# EXACT CONVERGENCE RATE AND LEADING TERM IN CENTRAL LIMIT THEOREM FOR STUDENT'S $T$ STATISTIC

### By Peter Hall and Qiying Wang

*Australian National University, Australian National University and University of Sydney*

The leading term in the normal approximation to the distribution of Student's $t$ statistic is derived in a general setting, with the sole assumption being that the sampled distribution is in the domain of attraction of a normal law. The form of the leading term is shown to have its origin in the way in which extreme data influence properties of the Studentized sum. The leading-term approximation is used to give the exact rate of convergence in the central limit theorem up to order $n^{-1/2}$, where $n$ denotes sample size. It is proved that the exact rate uniformly on the whole real line is identical to the exact rate on sets of just three points. Moreover, the exact rate is identical to that for the non-Studentized sum when the latter is normalized for scale using a truncated form of variance, but when the corresponding truncated centering constant is omitted. Examples of characterizations of convergence rates are also given. It is shown that, in some instances, their validity uniformly on the whole real line is equivalent to their validity on just two symmetric points.

**1. Introduction.** The Studentized mean is an early example of one of the most common approaches to adaptive statistical inference, where a nuisance parameter is replaced by its estimator and the effect on inference carefully gauged. Initially, in the case of Student's $t$ statistic, this was done under the assumption that the sampled distribution was normal, but later there developed a substantial literature, to which Gayen (1949, 1950, 1952) and Hyrenius (1950) were early contributors, on the effect of nonnormality on properties of the statistic. Wallace (1958), Bowman, Beauchamp and Shenton (1977) and Cressie (1980) have reviewed work in this area. Even in the case of normal data, where tables of the exact distribution have long been readily available, the issue of convergence (to normality) of the distribution of the $t$ statistic has been of both theoretical and practical interest for many years; see, for example, Anscombe (1950) and Gayen (1952).

From a theoretical viewpoint the problem of determining exact convergence rates for the $t$ statistic can be a particularly awkward one. Despite the statistic's simple representation in terms of the mean and mean of squares of independent data, its distribution is surprisingly difficult to approximate using methods for sums of independent random variables. The problem has, of course, long been solved under sufficiently severe moment conditions, but its treatment in more theoretically interesting cases, when its distribution is asymptotically normal but few other assumptions are made, is far from straightforward.

In a major advance, Bentkus and Götze (1996) gave bounds of general Berry–Esseen type for rates of convergence in the central limit theorem for Student's $t$ statistic when the data are independent and identically distributed. See also Chibisov (1980, 1984) and Slavova (1985). Bentkus, Bloznelis and Götze (1996) extended Bentkus and Götze's arguments to nonidentically distributed summands. Hall (1987) had earlier established Edgeworth expansions under moment conditions that were no more severe than existence of the moments actually appearing in the expansions. See also van Zwet (1984), Friedrich (1989), Putter and van Zwet (1998), Bentkus, Götze and van Zwet (1997), Wang and Jing (1999), Wang, Jing and Zhao (2000) and Bloznelis and Putter (1998, 2002).

However, moment conditions, even finite variance, are not the main prerequisite for convergence of the distribution of Student's $t$ statistic. In particular, Giné, Götze and Mason (1997) showed that a necessary and sufficient condition for the Studentized mean to have a limiting standard normal distribution is that the sampled distribution lie in the domain of attraction of the normal law. See also Logan, Mallows, Rice and Shepp (1973), Griffin and Mason (1991) and Egorov (1996). Although it is not of direct relevance to our work, we mention that the case where the data are from a time series is more complex. There, convergence in the conventional, deterministically normalized central limit theorem is not equivalent to convergence in the randomly normalized case; see Hahn and Zhang (1998).

In the present paper we assume no more than that the sampled distribution lies in the domain of attraction of the normal law, and describe rates of convergence, in the independent-data case, without reference to moment properties. We give the leading term in a normal approximation to the distribution of Student's $t$ statistic, and show that its form is strongly influenced by the effects that large data have on the statistic. Using the leading term, we derive the exact convergence rate in the central limit theorem, up to terms of order $n^{-1/2}$ (where $n$ denotes sample size), or up to order $n^{-1}$ when the sampled distribution satisfies Cramér's continuity condition.

We show that, if the third moment should happen to be finite, the leading term transforms into the conventional first term in an Edgeworth expansion of the distribution of Student's $t$ statistic. More generally, however, the leading term can be used to show that the exact rate of convergence over the whole real line is equivalent to the exact rate of convergence over very small sets, containing no more than three points. The number of points can be reduced to two if we seek necessary and sufficient characterizations of the convergence rate, rather than the exact rate itself. We draw connections to the rate of convergence of the distribution of a conventionally normalized, non-Studentized mean.

**2. Main results.** Let $X_1, X_2, \ldots$ be independent and identically distributed random variables, and let $X$ have the distribution of a generic $X_i$. Student's $t$ statistic, with numerator centered at its expectation, is defined to be

$$(2.1) \qquad T = \left( \sum_{i=1}^{n} X_i \right) \Big/ \left\{ \sum_{i=1}^{n} X_i^2 - n^{-1} \left( \sum_{i=1}^{n} X_i \right)^2 \right\}^{1/2}.$$

An alternative, more classical definition of the Studentized mean, in which the sample variance has divisor $n - 1$ rather than $n$, has the formula $(1 - n^{-1})^{-1/2} T$; see Gossett (1908). All our results hold for this version of Student's statistic, as well as that given by (2.1). The principal results are Theorems 2.1 and 2.2, which respectively describe the leading term and its role in a normal approximation to the distribution of $T$. Propositions 3.1 and 3.2 in the next section reveal the origins of the leading term, and in particular link it to the way in which extremes affect the distribution of $T$.

Write $\Phi$ and $\phi$ for the standard normal distribution and density functions, respectively. Put $b_n = \sup\{x : nx^{-2} E[X^2 I(|X| \leq x)] \geq 1\}$ and

$$(2.2) \qquad L_n(x) = nE(\Phi[x\{1 + (X/b_n)^2\}^{1/2} - (X/b_n)] - \Phi(x)).$$

THEOREM 2.1. *If the distribution of $X$ is in the domain of attraction of the normal law, and $E(X) = 0$, then*

$$(2.3) \qquad \sup_{-\infty < x < \infty} |P(T \leq x) - \{\Phi(x) + L_n(x)\}| = o(\delta_n) + O(n^{-1/2}).$$

*If, in addition, Cramér's condition holds, that is,*

$$\limsup_{|t| \to \infty} |E(e^{itX})| < 1,$$

*then $O(n^{-1/2})$ on the right-hand side of (2.3) may be replaced by $O(n^{-1})$.*

We noted in Section 1 that $T$ has a limiting standard normal distribution if and only if the distribution of $X$ is in the domain of attraction of the normal law and $E(X) = 0$. Theorem 2.1 argues that $L_n(x)$ is a leading term in an expansion of the distribution of $T$. As Theorem 2.2 will show, the exact order of magnitude of $L_n(x)$ is that of

$$(2.4) \qquad \begin{aligned} \delta_n = {} & nP(|X| > b_n) + nb_n^{-1}|E\{XI(|X| \leq b_n)\}| \\ & + nb_n^{-3}|E\{X^3I(|X| \leq b_n)\}| + nb_n^{-4}E\{X^4I(|X| \leq b_n)\}. \end{aligned}$$

THEOREM 2.2.   *Assume the distribution of $X$ is in the domain of attraction of the normal law and $E(X) = 0$. Then $\delta_n \to 0$ and*

$$(2.5) \qquad \sup_{-\infty < x < \infty} |L_n(x)| \asymp \delta_n$$

*as $n \to \infty$. Here and below, $a_n \asymp b_n$ denotes that*

$$0 < \liminf_{n \to \infty} a_n/b_n \leq \limsup_{n \to \infty} a_n/b_n < \infty.$$

*Property* (2.5) *continues to hold if the supremum over all $x$ is replaced by the supremum over $x \in \{-x_0, x_0, x_1\}$, where $x_0 > 3^{1/2}$ and $x_1$ is any real number not equal to $\pm x_0$. Furthermore, if $E(|X|^3) < \infty$, $E(X^2) = 1$ and $E(X^3) = \gamma$, then*

$$(2.6) \qquad \sup_{-\infty < x < \infty} |n^{1/2}L_n(x) - \tfrac{1}{6}\gamma(2x^2 + 1)\phi(x)| \to 0$$

*as $n \to \infty$.*

There exist examples of distributions in the domain of attraction of the normal law having zero mean and, for which any given one of the four components in the definition of $\delta_n$, at (2.4), dominate all the others along a subsequence. It follows that none of the terms of which $\delta_n$ is composed can be dropped if we require a full account of the rate of convergence in the central limit theorem. Formula (2.6) shows that in the case of finite third moment, the leading term is asymptotic to its conventional form in an Edgeworth expansion.

Together, properties (2.3) and (2.5) give concise results about the rate of convergence in the central limit theorem. For example, if $X$ is in the domain of attraction of the normal law, and $E(X) = 0$, then (2.3) and (2.5) imply that

$$(2.7) \qquad \sup_{-\infty < x < \infty} |P(T \leq x) - \Phi(x)| + n^{-1/2} \asymp \delta_n + n^{-1/2};$$

and $n^{-1/2}$ may be replaced by $n^{-1}$ if Cramér's condition is satisfied. One application to which (2.7) can be put is the derivation of characterizations of rates of convergence in the central limit theorem. In this regard, some examples can be found in Hall and Wang (2003), on which the present paper is based.

We conclude this section by mentioning that the convergence rate $\delta_n$ is the same as that in the case of the standard (i.e., non-Studentized) central limit theorem, where a sum of independent and identically distributed random variables is standardized for scale using $b_n$, but is centered conventionally, not using a truncated mean. That is, if we define $S_1 = b_n^{-1} \sum_{i \leq n} X_i$, $F_j(x) = P(S_j \leq x)$ and

$$(2.8) \qquad L_{n1}(x) = nE\{\Phi(x - X/b_n) - \Phi(x)\} - \tfrac{1}{2}nb_n^{-2}\phi'(x),$$

then, provided the distribution of $X$ is in the domain of attraction of the normal law and $E(X) = 0$, it is true that $\sup_{-\infty < x < \infty} |L_{n1}(x)| \asymp \delta_n$ and

$$(2.9) \qquad \sup_{-\infty < x < \infty} |F_1(x) - \{\Phi(x) + L_{n1}(x)\}| = o(\delta_n) + O(n^{-1/2}).$$

The methods of proof are similar to those given in Chapter 2 of Hall (1982). Alternatively, if we put $\sigma_n^2 = E\{X^2 I(|X| \leq b_n)\}$ and $S_2 = (\sum_{i \leq n} X_i)/(n^{1/2}\sigma_n)$, and define $L_{n2}(x)$ as at (2.8) but with $b_n$ there replaced by $n^{1/2}\sigma_n$, then (2.9) continues to hold if we replace $(F_1, L_{n1})$ by $(F_2, L_{n2})$.

The similarities between the Studentized and non-Studentized cases do not penetrate deeply, however. The leading terms in the respective settings are quite different. In the case of finite third moment, the leading terms are asymptotic to their respective Edgeworth forms, which are well known to have intrinsically different formulae.

## 3. Proofs.

3.1. *Proof of Theorem* 2.1. Let $\alpha > 0$ and define $Y_i = X_i I(|X_i| \leq \alpha b_n)$, $\rho_n = nP(|X| > \alpha b_n)$,

$$(3.1) \qquad \begin{aligned} \Psi_n(x) &= P\left[ \frac{\sum_{i \leq n} Y_i}{\{\sum_{i \leq n} Y_i^2 - n^{-1}(\sum_{i \leq n} Y_i)^2\}^{1/2}} \leq x \right], \\ M_{n1}(x) &= nE\{(\Phi[x\{1 + (X/b_n)^2\}^{1/2} - (X/b_n)] - \Phi(x))I(|X| > \alpha b_n)\}, \\ M_{n2}(x) &= nE\{(\Phi[x\{1 + (X/b_n)^2\}^{1/2} - (X/b_n)] - \Phi(x))I(|X| \leq \alpha b_n)\}. \end{aligned}$$

Theorem 2.1 is a direct consequence of the following two propositions, which will be proved in Sections 3.3 and 3.4.

PROPOSITION 3.1. *Assume the distribution of $X$ is in the domain of attraction of the normal law, and $E(X) = 0$. Then, for each $\alpha > 0$,*

$$(3.2) \qquad \sup_{-\infty < x < \infty} |P(T \leq x) - \{\Psi_n(x) + M_{n1}(x)\}| = o(\rho_n).$$

PROPOSITION 3.2. *Assume the distribution of $X$ is in the domain of attraction of the normal law, and $E(X) = 0$. Then, for each $\varepsilon > 0$ we have, for all sufficiently small $\alpha > 0$,*

$$(3.3) \qquad \sup_{-\infty < x < \infty} |\Psi_n(x) - \{\Phi(x) + M_{n2}(x)\}| \le \varepsilon \delta_n + O(n^{-1/2}).$$

*If, in addition, the distribution of $X$ satisfies Cramér's continuity condition, then the term $O(n^{-1/2})$ on the right-hand side of* (3.3) *may be replaced by $O(n^{-1})$.*

We remark that our method for proving Proposition 3.1 will show clearly that the leading-term fragment $M_{n1}$ derives principally from the largest summand among $X_1, \dots, X_n$, that is, from the value $X_{\max}$ of $X_i$ for which $|X_i|$ is greatest. Indeed, it may be proved that

$$M_{n1}(x) = E\{(\Phi[x\{1 + (X_{\max}/b_n)^2\}^{1/2} - (X_{\max}/b_n)] - \Phi(x))I(|X_{\max}| > \alpha b_n)\}$$
$$+ o(\rho_n),$$

uniformly in $x$. It follows that the leading term $L_n(x)$, introduced at (2.2) and defined as the limit of $M_{n1}$ as $\alpha \to 0$, also has this origin.

The connections to extremes arise in part through the major role that large summands play in convergence properties of series when the distribution of the summands has infinite variance. See Darling (1952), Arov and Bobrov (1960), Dwass (1966), Hall (1978), LePage, Woodroofe and Zinn (1981) and Resnick (1986) for discussion of more conventional settings. In the present case the main series where extremes cause difficulty is $\sum_{i \le n} X_i^2$, appearing in the definition of $T$ at (2.1). The summands here have finite variance if and only if the sampled distribution has finite fourth moment. However, extremes arising even from the series $\sum_{i \le n} X_i$ play a role in the leading term and so too in the convergence rate; see Hall (1984) for discussion of the latter issue.

3.2. *Proof of Theorem* 2.2. It is straightforward to show that $\delta_n \to 0$ and $\sup_{-\infty < x < \infty} |L_n(x)| = O(\delta_n)$. Therefore, it suffices to prove that

$$(3.4) \qquad \delta_n = O\left\{\sup_{x \in \mathcal{S}} |L_n(x)|\right\},$$

where $\mathcal{S} = \{-x_0, x_0, x_1\}$ is the set of three points in the statement of the theorem; and that (2.6) holds. This follows relatively straightforwardly.

3.3. *Proof of Proposition* 3.1. Put $V = \max_{i \le n} |X_i|$ and $J = \arg\max_{i \le n} |X_i|$; ties may be broken in any measurable way. Define $S$ to be the sign of $X_J$ and let $T_1 = \sum_{i \le n} X_i$, $T_2 = \sum_{i \le n} X_i^2$, $T_3 = \sum_{i \le n} Y_i + SVI(V > \alpha b_n)$ and $T_4 = \sum_{i \le n} Y_i^2 + V^2 I(V > \alpha b_n)$. The probability that two or more values of $|X_i|$, for $1 \le i \le n$, exceed $\alpha b_n$ equals $O(\rho_n^2)$. Therefore, $P\{(T_1, T_2) = (T_3, T_4)\} = 1 - O(\rho_n^2)$, whence it follows that, uniformly in $x$,

$$(3.5) \qquad \begin{aligned} P(T \le x) &= P\left\{ \frac{T_1}{(T_2 - n^{-1}T_1^2)^{1/2}} \le x \right\} \\ &= P\left\{ \frac{T_3}{(T_4 - n^{-1}T_3^2)^{1/2}} \le x \right\} + O(\rho_n^2). \end{aligned}$$

Put $\pi(v) = P(X \ge 0 | |X| = v)$. Conditional on $X_1, \ldots, X_n$, let $S(V)$ denote a random variable that takes the values $+1$ and $-1$ with probabilities $\pi(V)$ and $1 - \pi(V)$, respectively. Let $T_5 = \sum_{i \le n} Y_i + S(V)VI(V > \alpha b_n)$. Then $(T_3, T_4)$ has the same joint distribution as $(T_5, T_4)$, and so by (3.5) we have, uniformly in $x$,

$$(3.6) \qquad P(T \le x) = P(W \le x) + O(\rho_n^2),$$

where $W = T_5/(T_4 - n^{-1}T_5^2)^{1/2}$.

Define

$$\begin{aligned} T_Y &= \sum_{i \le n}(Y_i - EY_i), & T_B &= \sum_{i \le n}(Y_i^2 - EY_i^2), \\ \nu &= E\{XI(|X| \le \alpha b_n)\}, & \tau^2 &= E\{X^2 I(|X| \le \alpha b_n)\}. \end{aligned}$$

Note that a formula for $\Psi_n(x)$, equivalent to (3.1), is

$$(3.7) \qquad \Psi_n(x) = P\left[ \frac{T_Y + n\nu}{\{T_B + n\tau^2 - n^{-1}(T_Y + n\nu)^2\}^{1/2}} \le x \right].$$

Let the random variable $N_1$ have the standard normal distribution. The joint distribution of the vector $(b_n^{-1} T_Y, b_n^{-2} T_B)$, conditional on $V > \varepsilon b_n$, converges to the joint distribution of $(N_1, 0)$. In particular, the second component of the limiting distribution is degenerate at 0. The convergence has the following property: for all $\varepsilon > 0$,

$$(3.8) \qquad \sup_{v > \alpha b_n} \sup_{-\infty < x < \infty} |P(b_n^{-1} T_Y \le x; b_n^{-2}|T_B| \le \varepsilon | V = v) - P(N_1 \le x)| \to 0.$$

For a formal proof of (3.8), it suffices to observe that the joint distribution of $(\sum_{i \le n} Y_i, \sum_{i \le n} Y_i^2)$, conditional on $V = v > \alpha b_n$, equals the unconditional joint distribution of $(\sum_{i \le n-1} Y_i, \sum_{i \le n-1} Y_i^2)$; and that $b_n^{-1} \sum_{i \le n-1}(Y_i - EY_i) \to N_1$ in distribution, $b_n^{-2} \sum_{i \le n-1}(Y_i^2 - EY_i^2) \to 0$ in probability and $b_n^{-1}|E(Y_1)| + b_n^{-2} E(Y_1^2) \to 0$.

Since $T_4 = T_B + n\tau^2 + V^2 I(V > \alpha b_n)$, $T_5 = T_Y + n\nu + S(V)VI(V > \alpha b_n)$, $b_n^{-2} n\tau^2 \to 1$ and $b_n^{-1} n\nu \to 0$, then, for all $\varepsilon > 0$, we have from (3.8),

$$\sup_{v > \alpha b_n} \sup_{-\infty < x < \infty} |P[b_n^{-1}\{T_5 - S(V)(V/b_n)\} \le x;$$

$$|b_n^{-2} T_4 - 1 - (V/b_n)^2| \le \varepsilon | V = v] - P(N_1 \le x)| \to 0.$$

Therefore, if $N_2$ denotes a standard normal random variable that is independent of $V$, then

$$\sup_{v > \alpha b_n} \sup_{-\infty < x < \infty} \left| P(W \le x | V = v) - P\left[ \frac{N_2 + S(V)(V/b_n)}{\{1 + (V/b_n)^2\}^{1/2}} \le x \Big| V = v \right] \right| \to 0.$$

Equivalently,

$$P(W \le x | V = v) - \Phi[x\{1 + (v/b_n)^2\}^{1/2} - S(v)(v/b_n)] \to 0,$$

uniformly in $v > \alpha b_n$ and $-\infty < x < \infty$. Multiply throughout by $dF_n(v)$, where $F_n$ denotes the distribution function of $V$ conditional on $V > \alpha b_n$; integrate over $v > \alpha b_n$; and then multiply by $P(V > \alpha b_n)$, to prove that, uniformly in $x$,

$$
\begin{aligned}
(3.9) \quad & P(W \le x; V > \alpha b_n) \\
& = E(\Phi[x\{1 + (V/b_n)^2\}^{1/2} - S(V)(V/b_n)]I(V > \alpha b_n)) + o(\rho_n) \\
& = nE(\Phi[x\{1 + (X/b_n)^2\}^{1/2} - (X/b_n)]I(X > \alpha b_n)) + o(\rho_n).
\end{aligned}
$$

To derive the last identity, reformulate the expectation using an integration by parts argument, and note that

$$P\{S(V)V > y\} = nP(X > y) + O[\{nP(X > y)\}^2],$$

$$P\{S(V)V \le y\} = nP(X \le y) + O[\{nP(X \le y)\}^2],$$

where both remainders are of the stated orders uniformly in $y > \alpha b_n$ and $y < -\alpha b_n$, respectively.

Furthermore,

$$
\begin{aligned}
& E\{P(W \le x | V \le \alpha b_n)I(V \le \alpha b_n)\} \\
& \quad = E\left( I\left[ \frac{T_B + n\nu}{\{T_Y + n\tau^2 - n^{-1}(T_B + n\nu)^2\}^{1/2}} \le x \right] I(V \le \alpha b_n) \right) \\
& \quad = \Psi_n(x) - E\left( I\left[ \frac{T_B + n\nu}{\{T_Y + n\tau^2 - n^{-1}(T_B + n\nu)^2\}^{1/2}} \le x \right] I(V > \alpha b_n) \right),
\end{aligned}
$$

using (3.7) to obtain the last identity. A simpler version of the argument leading to (3.9) may be used to prove that the subtracted term above equals $\Phi(x)\rho_n + o(\rho_n)$, uniformly in $x$. Therefore,

$$(3.10) \qquad P(W \le x; V \le \alpha b_n) = \Psi_n(x) - \Phi(x)\rho_n + o(\rho_n),$$

uniformly in $x$. Combining (3.10) with (3.6) and (3.9), we conclude that (3.2) holds.

3.4. *Proof of Proposition* 3.2. Define $S_n = \sum_j X_j$, $S_n^* = \sum_j Y_j$, $V_n^2 = \sum_j X_j^2$ and $V_n^{*2} = \sum_j Y_j^2$. It is well known [e.g., Efron (1969)] that, for $x \geq 0$,

$$(3.11) \qquad \Psi_n(x) = P[S_n^*/V_n^* \leq x\{n/(n+x^2)\}^{1/2}].$$

Noting also that for $|u| \leq 1$,

$$\sup_{-\infty < x < \infty} |\Phi\{x(1+u^2)^{1/2} - u\}$$

$$- [\Phi(x) + \{-u + \tfrac{1}{6}u^3(2x^2+1) + \tfrac{1}{12}u^4 x(x^2-3)\}\phi(x)]| \leq C|u|^5,$$

where $C$ is an absolute constant, and that $nE\{|X/b_n|^5 I(|X| \leq \alpha b_n)\} \leq \alpha \delta_n$, we have that, for any $\alpha > 0$,

$$(3.12) \qquad \sup_{-\infty < x < \infty} |M_{n2}(x) - Q_{n1}(x)| \leq C\alpha \delta_n,$$

where $u_{nj} = nE\{(X/B_n)^j I_{(|X| \leq \alpha b_n)}\}$ and

$$Q_{n1}(x) = -u_{n1}\phi(x) + u_{n3}\tfrac{1}{6}(2x^2+1)\phi(x) + u_{n4}\tfrac{1}{12}x(x^2-3)\phi(x).$$

In view of (3.11) and (3.12), Proposition 3.2 will follow if, for each $\varepsilon > 0$, we have for all sufficiently small $\alpha > 0$,

$$(3.13) \qquad \sup_{-\infty < x < \infty} |P(S_n^*/V_n^* \leq x) - \{\Phi(x) + Q_{n1}(x)\}| \leq \varepsilon \delta_n + O(n^{-1/2}),$$

and $O(n^{-1/2})$ may be replaced by $O(n^{-1})$ if Cramér's condition is satisfied.

Without loss of generality, $x \geq 0$. Since the distribution of $X$ is in the domain of attraction of the normal law, then $\{S_n/V_n\}$ is stochastically bounded [see, e.g., Giné, Götze and Mason (1997)] and similarly $\{S_n^*/V_n^*\}$ is also stochastically bounded. Hence, by Theorem 2.5 of Giné, Götze and Mason (1997), for $x > \delta_n^{-1/12}$,

$$P(S_n^* \geq xV_n^*) \leq e^{-x}\sup_n E\{\exp(|S_n^*/V_n^*|)\} \leq A\exp(-\delta_n^{-1/12}) \leq A\delta_n^2.$$

(Here and below, $A$ denotes a positive constant which might be different at each appearance.) Moreover, $|1 - \Phi(x) - Q_{n1}(x)| \leq A\delta_n^2$ uniformly in $x > \delta_n^{-1/12}$. Therefore, (3.13) will follow if, for each $\varepsilon > 0$, we have for all sufficiently small $\alpha > 0$,

$$(3.14) \qquad \sup{}'|P(S_n^*/V_n^* \leq x) - \{\Phi(x) + Q_{n1}(x)\}| \leq \varepsilon \delta_n + O(n^{-1/2}),$$

and $O(n^{-1/2})$ may be replaced by $O(n^{-1})$ if Cramér's condition is satisfied, where $\sup'$ denotes the supremum over $x \in [0, \delta_n^{-1/12}]$.

Let $B_n^2 = nEY_1^2$ and $W_n = B_n^{-2}\sum_j(Y_j^2 - EY_j^2)$. Noting that $(1+y)^{1/2} = 1 + \tfrac{1}{2}y - \tfrac{1}{8}y^2 + \tfrac{1}{16}y^3 + \theta y^4$, where $\theta = \theta(y)$ satisfies $|\theta| \leq \tfrac{1}{16}$ for $|y| \leq \tfrac{1}{2}$, we

may prove that

$$P(S_n^*/V_n^* \le x)$$
$$= P\{S_n^* \le x B_n (1 + W_n)^{1/2}\}$$
$$\ge -P(|W_n| \ge \tfrac{1}{2}) + P\{S_n^* \le x B_n (1 + \tfrac{1}{2}W_n - \tfrac{1}{8}W_n^2 + \tfrac{1}{16}W_n^3 - \tfrac{1}{16}W_n^4)\},$$
$$P(S_n^*/V_n^* \le x)$$
$$= P\{S_n^* \le x B_n (1 + W_n)^{1/2}\}$$
$$\le P(|W_n| \ge \tfrac{1}{2}) + P\{S_n^* \le x B_n (1 + \tfrac{1}{2}W_n - \tfrac{1}{8}W_n^2 + \tfrac{1}{16}W_n^3 + \tfrac{1}{16}W_n^4)\}.$$

In view of Markov's inequality, it is readily seen that, for each $\varepsilon > 0$, we have for any sufficiently small $\alpha > 0$ and all sufficiently large $n$,

$$P(|W_n| \ge 1/2) \le 16 E(W_n^4) \le A(\alpha^4 \delta_n + \delta_n^2) \le \varepsilon \delta_n.$$

Hence, (3.14) will follow if we prove that, for each $\varepsilon > 0$, we have for $|\theta| \le 1/16$ and any sufficiently small $\alpha > 0$,

$$(3.15) \quad \begin{aligned} \sup{}'|P\{S_n^* \le x B_n (1 + \tfrac{1}{2}W_n - \tfrac{1}{8}W_n^2 + \tfrac{1}{16}W_n^3 + \theta W_n^4)\} \\ - \{\Phi(x) + Q_{n1}(x)\}| \\ \le \varepsilon \delta_n + O(n^{-1/2}), \end{aligned}$$

and $O(n^{-1/2})$ may be replaced by $O(n^{-1})$ if Cramér's condition is satisfied.

Let $\sum_{j \ne k}$, $\sum_{j \ne k \ne l}$ and $\sum_{j \ne k \ne l \ne m}$ denote summations over pairs, triples and quadruples, respectively, of distinct integers between 1 and $n$. Put $Z_j = Y_j^2 - E Y_j^2$. Simple calculations show that

$$B_n^6 W_n^3 = \sum_{j=1}^n Z_j^3 + 3 \sum_{j \ne k} Z_j (Z_k^2 - E Z_k^2) + \sum_{j \ne k \ne l} Z_j Z_k Z_l + W_{n1},$$

$$B_n^8 W_n^4 = \sum_{j=1}^n Z_j^4 + 4 \sum_{j \ne k} Z_j (Z_k^3 - E Z_k^3) + 12 \sum_{j \ne k} (Z_j^2 - E Z_j^2)(Z_k^2 - E Z_k^2) + W_{n2},$$

where $W_{n1} = 3(n-1)E(Z_1^2)\sum_j Z_j$ and

$$\begin{aligned} W_{n2} = {} & 4(n-1)E(Z_1^3)\sum_{j=1}^n Z_j + 24(n-1)E(Z_1^2)\sum_{j=1}^n (Z_j^2 - E Z_j^2) \\ & + 12n(n-1)(E Z_1^2)^2 + 24 \sum_{j \ne k \ne l} Z_j Z_k Z_l^2 + \sum_{j \ne k \ne l \ne m} Z_j Z_k Z_l Z_m. \end{aligned}$$

Therefore,

$$P\left\{S_n^* \le x B_n \left(1 + \frac{1}{2}W_n - \frac{1}{8}W_n^2 + \frac{1}{16}W_n^3 + \theta W_n^4\right)\right\}$$

$$= P\left\{\frac{1}{B_n}\sum_{j=1}^n \xi_j(x) + \frac{x}{B_n^4}\sum_{j\neq k}\varphi_{jk} + \frac{x}{B_n^6}\sum_{j\neq k\neq l}\psi_{jkl}\right.$$

$$\left.\leq x(1 + W_{n3}) - \frac{nE\eta_1(x)}{B_n}\right\},$$

where $\xi_j(x) = \eta_j(x) - E\eta_j(x)$, $\psi_{jkl} = -\frac{1}{16}Z_j Z_k Z_l$, $W_{n3} = \frac{1}{16}B_n^{-6}W_{n1} + \theta B_n^{-8}W_{n2}$,

$$\eta_j(x) = Y_j - \frac{x}{2B_n}Z_j + \frac{x}{8B_n^3}Z_j^2 - \frac{x}{16B_n^5}Z_j^3 - \frac{\theta x}{B_n^7}Z_j^4,$$

$$\varphi_{jk} = \frac{1}{8}Z_j Z_k - \frac{3}{16B_n^2}Z_j(Z_k^2 - EZ_k^2) - \frac{4\theta}{B_n^4}Z_j(Z_k^3 - EZ_k^3)$$

$$- \frac{12\theta}{B_n^4}(Z_j^2 - EZ_j^2)(Z_k^2 - EZ_k^2).$$

It is readily seen that

$$E(B_n^{-6}W_{n1})^4 \leq A\delta_n^4(\alpha^4\delta_n + \delta_n^2) \quad \text{and} \quad E(B_n^{-8}W_{n2})^2 \leq A\delta_n^2(\alpha^2\delta_n + \delta_n^2).$$

Hence, for each $\varepsilon > 0$, we have for any sufficiently small $\alpha > 0$ and all sufficiently large $n$,

$$P(|W_{n3}| \geq 2\varepsilon\delta_n) \leq P(|B_n^{-6}W_{n1}| \geq \varepsilon\delta_n) + P(|B_n^{-8}W_{n2}| \geq \varepsilon\delta_n)$$

$$\leq A\{\varepsilon^{-4}(\alpha^4\delta_n + \delta_n^2) + \varepsilon^{-2}(\alpha^2\delta_n + \delta_n^2)\} \leq \varepsilon\delta_n.$$

Result (3.15) now follows easily from the following three propositions. We will only prove Propositions 3.3 and 3.4 in subsequent sections. The proof of Proposition 3.5 is relatively straightforward although requiring tedious algebra, and hence details are omitted. The proof of Proposition 3.2 is therefore complete.

PROPOSITION 3.3. *For all* $0 < \alpha \leq \frac{1}{2}$,

(3.16)
$$\sup{}' \sup_{-\infty < y < \infty} \left| P\left\{\frac{1}{B_n}\sum_{j=1}^n \xi_j(x) + \frac{x}{B_n^4}\sum_{j\neq k}\varphi_{jk} + \frac{x}{B_n^6}\sum_{j\neq k\neq l}\psi_{jkl} \leq y\right\}\right.$$

$$\left. - \{\Phi(y) + \mathcal{L}_n(y)\}\right|$$

$$= o(\delta_n) + O(n^{-1/2}),$$

*where* $\mathcal{L}_n(y) = n[E\Phi\{y - \xi_1(x)/B_n\} - \Phi(y)] - \frac{1}{2}\Phi^{(2)}(y)$.

PROPOSITION 3.4. *If* $\limsup_{|t|\to\infty} |Ee^{itX}| < 1$, *then the term* $O(n^{-1/2})$ *on the right-hand side of* (3.16) *may be replaced by* $O(n^{-1})$.

PROPOSITION 3.5. *For each $\varepsilon > 0$, we have for any sufficiently small $\alpha > 0$,*

$$(3.17) \quad \sup{}'|\Phi[x - \{nE\eta_1(x)/B_n\}] - \Phi(x) + Q_{n2}(x)| \leq \varepsilon \delta_n + O(n^{-1}),$$

$$(3.18) \quad \sup{}'|\mathcal{L}_n[x - \{nE\eta_1(x)/B_n\}] - Q_{n3}(x)| \leq \varepsilon \delta_n + O(n^{-1}),$$

*where $u_{nj} = nE\{(X/B_n)^j I(|X| \leq \alpha b_n)\}$, $Q_{n2}(x) = u_{n1}\phi(x) + \frac{1}{8}xu_{n4}\phi(x)$ and*

$$Q_{n3}(x) = u_{n3}\tfrac{1}{6}(2x^2 + 1)\phi(x) + u_{n4}\tfrac{1}{24}x(2x^2 - 3)\phi(x).$$

3.5. *Proof of Proposition* 3.3. Standard methods based on Taylor's expansion, although requiring tedious algebra, may be used to establish the following lemmas. Define

$$u_{nj} = nE\{(X/B_n)^j I(|X| \leq \alpha b_n)\}, \qquad g(t,x) = E[\exp\{it\xi_1(x)/B_n\}]$$

and

$$f_n(t,x) = e^{-t^2/2}[1 + n\{g(t,x) - 1\} + \tfrac{1}{2}t^2].$$

LEMMA 3.6. *If $0 < \alpha \leq \frac{1}{2}$, then for all sufficiently large $n$,*

$$(3.19) \quad |nB_n^{-2}E\xi_1^2(x) - (1 - xu_{n3} + \tfrac{1}{4}x^2u_{n4})| \leq 2(1 + x^2)(\alpha\delta_n + n^{-1}),$$

$$(3.20) \quad |nB_n^{-3}E\xi_1^3(x) - (u_{n3} - \tfrac{3}{2}xu_{n4})| \leq 12(1 + |x|^3)(\alpha\delta_n + n^{-1}),$$

$$(3.21) \quad |nB_n^{-4}E\xi_1^4(x) - u_{n4}| \leq 32(1 + x^4)(\alpha\delta_n + n^{-1}),$$

$$(3.22) \quad nB_n^{-5}E|\xi_1(x)|^5 \leq 32(1 + |x|^5)\alpha\delta_n.$$

LEMMA 3.7. *There exists a constant $c_0 > 0$ such that, for all $\alpha \in (0, \frac{1}{2}]$, $|t| \leq c_0 n^{1/2}$, all $x \in [0, \delta_n^{-1/12}]$ and all sufficiently large $n$,*

$$(3.23) \quad |g(t,x)| \leq e^{-t^2/8n},$$

$$(3.24) \quad |g^n(t,x) - e^{-t^2/2}| \leq A(1 + x^4)(1 + \alpha^{-1})\delta_n(t^2 + t^4)e^{-t^2/8},$$

$$(3.25) \quad |g^n(t,x) - f_n(t,x)| \leq \{A(1 + x^8)(1 + \alpha^{-2})\delta_n^2(t^4 + t^8) + 2n^{-1}t^4\}e^{-t^2/8}.$$

Throughout the proof of Proposition 3.3, we assume that $0 < \alpha \leq \frac{1}{2}$, $0 \leq x \leq \delta_n^{-1/12}$ and $n$ is sufficiently large. Define $\bar{\varphi}_{jk} = \varphi_{jk} + \varphi_{kj}$, $T_n = B_n^{-1}\sum_j \xi_j(x)$ and

$$\Delta_{n,m} = \frac{x}{B_n^4}\sum_{j=1}^{m-1}\sum_{k=j+1}^{n}\bar{\varphi}_{jk} + \frac{6x}{B_n^6}\sum_{j=1}^{m-2}\sum_{k=j+1}^{n}\sum_{l=k+1}^{n}\psi_{jkl}.$$

Noting that $B_n^2 = nEY_1^2 = b_n^2$ for sufficiently large $n$, we obtain that $|Y_j| \leq B_n/2$ and $|Z_j| \leq B_n^2/2$. Using these properties, $E(\psi_{123}|X_j) = 0, j = 1, 2, 3,$

and $E(\bar{\varphi}_{12}|X_j) = 0, j = 1, 2$, we may deduce that $E(\varphi_{12}^2) \leq 2^8 (EY_1^4)^2$, $E(\psi_{123}^2) \leq (EY_1^4)^3$ and, for $1 \leq m \leq n$,

$$
\begin{aligned}
(3.26) \quad E(\Delta_{n,m}^2) &\leq 2x^2 \{ mn B_n^{-8} E(\varphi_{12}^2) + mn^2 B_n^{-12} E(\psi_{123}^2) \} \\
&\leq 2^{10} mn^{-1} x^2 \delta_n^2,
\end{aligned}
$$

the last inequality following from the fact that $nB_n^{-4} EY_1^4 \leq \delta_n$. Moreover, noting that $n^{-1} \leq \delta_n \to 0$, $|u_{n4}| \leq \delta_n$ and $|u_{n3}| \leq (1 + \alpha^{-1})\delta_n$, it follows easily from (3.19)–(3.21) that

$$
(3.27) \quad \left| \frac{E\xi_1^2(x)}{EY_1^2} - 1 \right| \leq 2(1 + x^2)(1 + \alpha^{-1})\delta_n \leq \frac{1}{2},
$$

$$
(3.28) \quad nB_n^{-3} |E\xi_1^3(x)| \leq 32(1 + |x|^3)(1 + \alpha^{-1})\delta_n,
$$

$$
(3.29) \quad nB_n^{-4} E\xi_1^4(x) \leq 65(1 + x^4)\delta_n.
$$

We now turn back to the proof of (3.16). Using the identities

$$
\Delta_{n,n} = \frac{x}{B_n^4} \sum_{j \neq k} \varphi_{jk} + \frac{x}{B_n^6} \sum_{j \neq k \neq l} \psi_{jkl}
$$

and $\int e^{ity} d\{\Phi(y) + \mathcal{L}_n(y)\} = f_n(t, x)$, and Esseen's smoothing lemma [e.g., Petrov (1975), page 109], it may be shown that

$$
\begin{aligned}
(3.30) \quad &\sup_{-\infty < y < \infty} |P(T_n + \Delta_{n,n} \leq y) - \{\Phi(y) + \mathcal{L}_n(y)\}| \\
&\leq \int_{|t| \leq \min\{\delta_n^{-2}, c_0 n^{1/2}\}} |E \exp\{it(T_n + \Delta_{n,n})\} - f_n(t,x)| |t|^{-1} \, dt \\
&\quad + A(\delta_n^2 + n^{-1/2}) \sup_{-\infty < y < \infty} |(d/dy)\{\Phi(y) + \mathcal{L}_n(y)\}| \\
&\leq \sum_{j=1}^{4} I_{jn} + A(\delta_n^2 + n^{-1/2})(1 + \alpha^{-1}\delta_n^{2/3}),
\end{aligned}
$$

where $c_0$ is as in Lemma 3.7,

$$
\begin{aligned}
I_{1n} = \int_{|t| \leq \delta_n^{-1/4}} &|E \exp\{it(T_n + \Delta_{n,n})\} \\
&- E \exp(itT_n) - itE\{\Delta_{n,n} \exp(itT_n)\}| |t|^{-1} \, dt,
\end{aligned}
$$

$$
\begin{aligned}
I_{2n} = \int_{|t| \leq \delta_n^{-1/4}} &2|E \exp(itT_n) - f_n(t,x)| |t|^{-1} \, dt \\
&+ \int_{\delta_n^{-1/4} \leq |t| \leq c_0 n^{1/2}} |E \exp(itT_n)| |t|^{-1} \, dt,
\end{aligned}
$$

$$
I_{3n} = \int_{|t| \leq \delta_n^{-1/4}} |E\{\Delta_{n,n} \exp(itT_n)\}| \, dt,
$$

$$
I_{4n} = \int_{\delta_n^{-1/4} \leq |t| \leq \min\{\delta_n^{-2}, c_0 n^{1/2}\}} |E \exp\{it(T_n + \Delta_{n,n})\}| |t|^{-1} \, dt,
$$

and we have used the property, implied by (3.27)–(3.29), that

$$\sup_{-\infty < y < \infty} |\mathcal{L}'_n(y)| \leq A\left\{ \frac{1}{2}\left| \frac{nE\xi_1^2(x)}{B_n^2} - 1 \right| + \frac{n|E\xi_1^3(x)|}{6B_n^3} + \frac{nE\xi_1^4(x)}{24B_n^4} \right\}$$

$$\leq A(1+x^4)(1+\alpha^{-1})\delta_n \leq A(1+\alpha^{-1})\delta_n^{2/3}.$$

Using (3.26) and the fact that $|e^{iu} - 1 - iu| \leq u^2/2$, it can be shown that

$$(3.31) \qquad I_{1n} \leq \tfrac{1}{2} \int_{|t| \leq \delta_n^{-1/4}} E(\Delta_{n,n}^2)|t|\, dt \leq 2^9 x^2 \delta_n^{3/2} \leq 2^9 \delta_n^{4/3}.$$

Using Lemma 3.7, we obtain

$$(3.32) \qquad I_{2n} \leq A\{(1+x^8)(1+\alpha^{-2})\delta_n^2 + n^{-1}\} \leq A\{(1+\alpha^{-2})\delta_n^{4/3} + n^{-1}\}.$$

Next we estimate $I_{3n}$ and $I_{4n}$. Treating the former first, note that $E(\varphi_{12}|X_1) = E(\varphi_{12}|X_2) = 0$ and $E\varphi_{12}^2 \leq 2^8(EY_1^4)^2$, and that as in Bickel, Götze and van Zwet (1986),

$$(3.33) \quad E(\varphi_{12}\exp[it\{\xi_1(x) + \xi_2(x)\}/B_n]) = -\frac{t^2}{B_n^2}E\{\xi_1(x)\xi_2(x)\varphi_{12}\} + l_1(x),$$

where by using $|e^{iu} - 1 - iu| \leq u^2/2$, $|e^{iu} - 1| \leq |u|$, (3.27) and (3.29),

$$|l_1(x)| \leq \left| E\left[ \varphi_{12}\left\{ e^{it\xi_1(x)/B_n} - 1 - \frac{it\xi_1(x)}{B_n} \right\}\left\{ e^{it\xi_2(x)/B_n} - 1 \right\} \right] \right|$$

$$+ \frac{|t|}{B_n}\left| E\left[ \varphi_{12}\xi_1(x)\left\{ e^{it\xi_2(x)/B_n} - 1 - \frac{it\xi_2(x)}{B_n} \right\} \right] \right|$$

$$\leq \frac{|t|^3}{2B_n^3}E[|\varphi_{12}|\{|\xi_1(x)|\xi_2^2(x) + \xi_1^2(x)|\xi_2(x)|\}]$$

$$\leq \frac{|t|^3}{B_n^3}(E\varphi_{12}^2)^{1/2}\{E\xi_1^2(x)\}^{1/2}\{E\xi_1^4(x)\}^{1/2}$$

$$\leq A(1+x^2)|t|^3 n^{-1/2}B_n^{-1}(EY_1^4)(EY_1^2)^{1/2}\delta_n^{1/2}$$

$$\leq A(1+x^2)|t|^3 n^{-2}\delta_n^{3/2}B_n^4,$$

since $nB_n^{-4}EY_1^4 \leq \delta_n$ and $B_n^2 = nEY_1^2$. Tedious but elementary calculation shows that

$$|E\{\xi_1(x)\xi_2(x)\varphi_{12}\}| \leq A(1+x^2)n^{-2}\delta_n^2 B_n^6.$$

Substituting into (3.33), we deduce that

$$(3.34) \quad |E(\varphi_{12}\exp[it\{\xi_1(x) + \xi_2(x)\}/B_n])| \leq A(1+x^2)(t^2 + |t|^3)n^{-2}\delta_n^{3/2}B_n^4.$$

Similarly, it follows from the identities $E(\psi_{123}|X_j) = 0$, for $j = 1, 2, 3$, and from $E\psi_{123}^2 \leq (EY_1^4)^3$ and (3.27), that

$$(3.35) \quad |E(\psi_{123}\exp[it\{\xi_1(x) + \xi_2(x) + \xi_3(x)\}/B_n])| \leq A|t|^3 n^{-3}\delta_n^{3/2}B_n^6.$$

From (3.34), (3.35) and (3.23), it can be seen that

$$|E\{\Delta_{n,n}\exp(itT_n)\}|$$
$$\leq |x|n^2 B_n^{-4}|E\{\varphi_{12}\exp(itT_n)\}| + |x|n^3 B_n^{-6}|E\{\psi_{123}\exp(itT_n)\}|$$
$$\leq A(1+|x|^3)\delta_n^{3/2}(t^2+|t|^3)e^{-t^2/8},$$

and hence

$$(3.36)\quad I_{3n} = \int_{|t|\leq\delta_n^{-1/4}}|E\{\Delta_{n,n}\exp(itT_n)\}|\,dt \leq A(1+|x|^3)\delta_n^{3/2} \leq A\delta_n^{5/4}.$$

We next estimate $I_{4n}$. Put $\Delta_{n,m}^* = \Delta_{n,n} - \Delta_{n,m}$. In view of (3.26),

$$|E\exp\{it(T_n+\Delta_{n,n})\} - E\exp\{it(T_n+\Delta_{n,m}^*)\} - itE\Delta_{n,m}\exp\{it(T_n+\Delta_{n,m}^*)\}|$$
$$\leq 2^9 t^2 x^2 mn^{-1}\delta_n^2.$$

This inequality, together with the independence of the $X_k$'s, implies that for any $1\leq m\leq n$,

$$(3.37)\quad \begin{aligned}&|E\exp\{it(T_n+\Delta_{n,n})\}|\\ &\qquad \leq |g(t,x)|^{m-2} + A|x|\delta_n|t||g(t,x)|^{m-5} + At^2x^2mn^{-1}\delta_n^2,\end{aligned}$$

where we have used the bound $E|\Delta_{n,m}| \leq (E|\Delta_{n,m}|^2)^{1/2} \leq A|x|\delta_n$.

Let $n_0 = [16nt^{-2}\log(\delta_n^{-1})] + 5$, where $[\cdot]$ denotes the integer part function. It is clear that $1 \leq n_0 \leq n$ for $\delta_n^{-1/4} \leq |t| \leq \min\{\delta_n^{-2}, c_0 n^{1/2}\}$, for $n$ large enough. Hence, choosing $m = n_0$ in (3.37) and using (3.23), we get

$$(3.38)\quad \begin{aligned}I_{4n} &= \int_{\delta_n^{-1/4}\leq|t|\leq\min\{\delta_n^{-2},c_0 n^{1/2}\}}|E\exp\{it(T_n+\Delta_{n,n})\}||t|^{-1}\,dt\\ &\leq A(1+x^2)(\log\delta_n^{-1})^2\delta_n^2 \leq A\delta_n^{4/3}.\end{aligned}$$

Substituting the bounds for $I_{1n},\ldots,I_{4n}$ back into (3.30), and recalling that $\delta_n \to 0$, we obtain (3.16), and hence complete the proof of Proposition 3.3.

3.6. *Proof of Proposition* 3.4. Without loss of generality we assume that $\delta_n \leq n^{-1/3}$. Indeed, for $\delta_n \geq n^{-1/3}$, it is obvious that the term $O(n^{-1/2})$ on the right-hand side of (2.3) can be replaced by $O(n^{-1})$. Note that $\delta_n \leq n^{-1/3}$ implies that $nP(|X| \geq b_n) \leq n^{-1/3}$. This, together with the fact that the distribution of $X$ is in the domain of attraction of a normal law, implies that $EX^2 < \infty$.

We continue to use the notation in the proof of Proposition 3.3. Further, we put

$$\tilde{\Delta}_{n,m}^{(1)} = \frac{x}{B_n^4}\sum_{j=1}^{m}\sum_{k=m+1}^{n}\bar{\varphi}_{jk} + \frac{6x}{B_n^6}\sum_{j=1}^{m}\sum_{k=m+1}^{n}\sum_{l=k+1}^{n}\psi_{jkl},$$

$$\tilde{\Delta}_{n,m}^{(2)} = \frac{x}{B_n^4}\sum_{j=m+1}^{n-1}\sum_{k=j+1}^{n}\bar{\varphi}_{jk} + \frac{6x}{B_n^6}\sum_{j=m+1}^{n-2}\sum_{k=j+1}^{n}\sum_{l=k+1}^{n}\psi_{jkl}.$$

As in the proof of (3.26), we have that, for $0 \leq x \leq \delta_n^{-1/12}$,

$$(3.39) \quad E(\Delta_{n,n} - \tilde{\Delta}_{n,m}^{(1)} - \tilde{\Delta}_{n,m}^{(2)})^2 \leq 2^{10} m^2 n^{-2} x^2 \delta_n^2 \leq Am^2 n^{-2} \delta_n^{11/6}.$$

Hence, for $m_0 = C \log n$, where $C$ is a constant that we shall specify later,

$$P(|\Delta_{n,n} - \tilde{\Delta}_{n,m_0}^{(1)} - \tilde{\Delta}_{n,m_0}^{(2)}| \geq n^{-1}) \leq AC^2 (\log n)^2 \delta_n^{11/6} \leq AC^2 \delta n^{3/2}.$$

Proposition 3.4 will now follow if we show that, for $0 \leq x \leq \delta_n^{-1/12}$,

$$(3.40) \quad \begin{aligned} &\sup_{-\infty < y < \infty} |P(T_n + \tilde{\Delta}_{n,m_0}^{(1)} + \tilde{\Delta}_{n,m_0}^{(2)} \leq y) - \{\Phi(y) + \mathcal{L}_n(y)\}| \\ &= o(\delta_n) + O(n^{-1}). \end{aligned}$$

Throughout the proof of Proposition 3.4, we assume that $0 < \alpha \leq \frac{1}{2}$, $0 \leq x \leq \delta_n^{-1/12}$ and $n$ is sufficiently large. We need the following lemma, the proof of which can be found in Prawitz ([1972]). See also Bentkus, Götze and van Zwet ([1997]).

LEMMA 3.8.   *Let $F$ be a distribution function with characteristic function $f$. Then for all $y \in R$ and $T > 0$, it holds that*

$$(3.41) \quad \lim_{z \downarrow y} F(z) \leq \tfrac{1}{2} + \text{P.V.} \int_{-T}^{T} \exp(-iyt) T^{-1} K(t/T) f(t) \, dt,$$

$$(3.42) \quad \lim_{z \uparrow y} F(z) \geq \tfrac{1}{2} - \text{P.V.} \int_{-T}^{T} \exp(-iyt) T^{-1} K(-t/T) f(t) \, dt,$$

*where*

$$\text{P.V.} \int_{-T}^{T} = \lim_{h \downarrow 0} \left( \int_{-T}^{-h} + \int_{h}^{T} \right),$$

*and $2K(s) = K_1(s) + iK_2(s)/(\pi s)$,*

$$K_1(s) = 1 - |s|, \qquad K_2(s) = \pi s(1 - |s|) \cot \pi s + |s| \qquad \text{for } |s| < 1,$$

*and $K(s) \equiv 0$ for $|s| \geq 1$.*

We shall give the proof of (3.40) by using Lemma 3.8 and some of the techniques of Bentkus, Götze and van Zwet ([1997]). By $E_k(\cdot) = E(\cdot|X_{k+1}, \ldots, X_n)$ we shall denote expectation conditional on $X_{k+1}, \ldots, X_n$. Define

$$\tau_1 = n^{1/2} \delta_n^{-2/3} E_{m_0} \left| \frac{x}{B_n^4} \sum_{k=m_0+1}^{n/2} \bar{\varphi}_{1k} \right|, \qquad \tau_2 = n^{1/2} \delta_n^{-2/3} E_{m_0} \left| \frac{x}{B_n^4} \sum_{k=n/2+1}^{n} \bar{\varphi}_{1k} \right|,$$

and put $\tau_0 = 1 - \limsup_{|t| \to \infty} |E e^{itX}|$,

$$H = \frac{n^{1/2} \delta_n^{-2/3} \tau_0}{16(1 + \tau_1 + \tau_2)}.$$

As in the proof of (3.26),

$$E(\tau_1 + \tau_2)^2 \le 2n\delta_n^{-4/3} \left\{ E\left( \frac{x}{B_n^4} \sum_{k=m_0+1}^{n/2} \bar{\varphi}_{1k} \right)^2 + E\left( \frac{x}{B_n^4} \sum_{k=n/2+1}^{n} \bar{\varphi}_{1k} \right)^2 \right\}$$

$$\le Ax^2 \delta_n^{2/3}.$$

This, together with the bound $0 \le x \le \delta_n^{-1/12}$, implies that

$$(3.43) \qquad E(H^{-2}) \le An^{-1} \delta_n^{4/3} E(1 + \tau_1 + \tau_2)^2 \le An^{-1} \delta_n^{4/3}.$$

Also, we have that $H \le (\tau_0/16) n^{1/2} \delta_n^{-2/3}$.

Returning to the proof of (3.40), note that $H$ depends only on $X_{m_0+1}, \ldots, X_n$. Using (3.41), and arguing as in Bentkus, Götze and van Zwet (1997), we obtain

$$(3.44) \qquad 2P(T_n + \tilde{\Delta}_{n,m_0}^{(1)} + \tilde{\Delta}_{n,m_0}^{(2)} \le y) \le 1 + EI_1 + EI_2,$$

where, with $f(t) = E_{m_0} \exp[it\{T_n + \tilde{\Delta}_{n,m_0}^{(1)} + \tilde{\Delta}_{n,m_0}^{(2)}\}]$,

$$I_1 = H^{-1} \int_R \exp(-iyt) K_1(t/H) f(t) \, dt,$$

$$I_2 = \frac{i}{\pi} \text{P.V.} \int_R \exp(-iyt) K_2(t/H) f(t) t^{-1} \, dt.$$

The following results are derived by Hall and Wang (2003), on which the present paper is based:

$$(3.45) \qquad\qquad\qquad |EI_1| = o(\delta_n) + O(n^{-1}),$$

$$(3.46) \qquad |EI_2 + 1 - 2\{\Phi(y) + \mathcal{L}_n(y)\}| = o(\delta_n) + O(n^{-1}).$$

It follows from (3.44)–(3.46) that

$$P(T_n + \tilde{\Delta}_{n,m_0}^{(1)} + \tilde{\Delta}_{n,m_0}^{(2)} \le y) \le \Phi(y) + \mathcal{L}_n(y) + o(\delta_n) + O(n^{-1}).$$

Similarly, using (3.42) and symmetry arguments, one can show that

$$P(T_n + \tilde{\Delta}_{n,m_0}^{(1)} + \tilde{\Delta}_{n,m_0}^{(2)} > y) \le 1 - \{\Phi(y) + \mathcal{L}_n(y)\} + o(\delta_n) + O(n^{-1}).$$

Result (3.40) now follows, and hence the proof of Proposition 3.4 is complete.

**Acknowledgment.**  The authors thank a referee and Editors for their valuable comments and suggestions.

Centre for Mathematics
and its Applications
Mathematical Sciences Institute
Australian National University
Canberra ACT 0200
Australia
e-mail: Peter.Hall@maths.anu.edu.au

Centre for Mathematics
and its Applications
Mathematical Sciences Institute
Australian National University
Canberra ACT 0200
Australia
and
School of Mathematics and Statistics
University of Sydney
Sydney NSW 2006
Australia
e-mail: Qiying.Wang@maths.anu.edu.au